\documentclass[12pt,a4paper]{article}
\usepackage[T2A]{fontenc}
\usepackage[cp1251]{inputenc}
\usepackage[english]{babel}
\usepackage{amssymb, amsmath, latexsym,amsthm}
\usepackage{amsfonts,mathtext,cite,enumerate,float}

\usepackage{authblk}

\usepackage{geometry} 
\title{Gorin's problem for individual simple partial fractions}

\author{Petr Chunaev\footnote {chunaev@itmo.ru; National Center for Cognitive Research, ITMO University} \, and Vladimir Danchenko\footnote{vdanch2012@yandex.ru; Vladimir State University named after Alexander and Nikolay Stoletovs}}

\date{}

\begin{document}

\maketitle

\begin{abstract}
The main result of the paper is a lower estimate for the moduli of imaginary parts of the poles of a simple partial fraction (i.e. the logarithmic derivative of an algebraic polynomial) under the condition that the $L^\infty(\mathbb{R})$-norm of the fraction is unit (Gorin's problem). In contrast to the preceding results, the estimate takes into account the residues associated with the poles.

Moreover, a new estimate for the moduli is obtained in the case when the  $L^\infty(\mathbb{R})$-norm of the derivative of the simple partial fraction is unit (Gelfond's problem).

\end{abstract}

\section{Introduction}
Gorin's problem is formulated as follows. Find a lower estimate for
$$
d_n(\mathbb R,p)=\inf\left\{Y(\rho_n):\;
\|\rho_n\|_{L^{p}({\mathbb R})}\le 1 \right\},\quad \quad
1<p\le\infty,
 \eqno{(1)}
$$
where $Y(\rho_n):=\min_{k=1,\ldots,m} |{\rm Im}\,{\xi}_k|$, and
$$
\rho_n(z):=\left(\ln
\prod_{k=1}^m(z-{\xi}_k)^{n_k}\right)'=\sum_{k=1}^m\frac{n_k}{z-{\xi}_k}
\eqno{(2)}
$$
is a {\it simple partial fraction} (SPF), i.e. the logarithmic derivative of an algebraic polynomial $Q(z)=\prod_{k=1}^m(z-{\xi}_k)^{n_k}$ of a given degree $n=\sum _{k=1}^m{n_k}$ (the number $n$ is called the {\it order} of the SPF). Note that SPFs, being discrete Cauchy potentials, have wide applications in Electrostatic Field Theory, Potential Theory (see \cite[\S 1]{[7]}, \cite[\S 3]{Ei}) and Theory of Partial Differential Equations \cite{[1]}.

We start with a brief history of the problem~(1). For $p=\infty$ (this is the most difficult case) the problem was considered in \cite{[1],[2],[3],[4],[5],[6]}. The question about the principal possibility of the estimate $d_n(\mathbb R,\infty)\ge c(n)>0$ was stated and positively resolved by Gorin~\cite{[1]}. Later, Nikolaev \cite{[2]} obtained the estimate
$$
d_n(\mathbb R,\infty)\ge 2\,(\sqrt{2}-1)^{n-1},\qquad n\in
{\mathbb N}.
$$
In \cite{[2]}, the following problem was also stated: \textit{Is it true that $d_n(\mathbb R,\infty)\to 0$ as $n\to\infty$?} An essential improvement for Nikolaev's estimate was done by Gelfond~\cite{[3]}:
$$
d_n(\mathbb R,\infty)\ge (17\ln n)^{-1},\qquad n\ge n_{0}.
$$
Furthermore, Katsnelson~\cite{[4]} obtained a certain improvement of this estimate but with the same logarithmic rate of minorant decrease. However, Nikolaev's problem remained open. The final solution on the whole class of SPFs~(2) was given\footnote{The expression $\alpha\asymp \beta$ means that there exist absolute positive constants $c$ and $C$ such that $c\beta \le \alpha\le C\beta$.} in \cite{[6]}:
$$
d_n(\mathbb R,\infty)\asymp \frac{\ln\ln n}{\ln n}.
\eqno{(3)}
$$
It was also shown there that for finite $p$ the value $d_n(\mathbb R,p)$ does not tend to zero and is bounded from below by a positive constant depending only on $p$. Actually the following stronger result holds \cite{[6]}:
$$
\| \rho^{\pm}_n \|_{L^{\infty}({\mathbb R})} \le \beta_p \cdot
 \|\rho_n\|_{L^{p}({\mathbb R})}^q,
 \qquad\hbox{where} \qquad
 \beta_p\le 2p\,\sin^{-q}({\pi}/{p}),\quad p^{-1}+q^{-1}=1,
$$
and $\rho^{\pm}_n$ denote partial sums from~(2) containing all poles from the half-planes ${\mathbb C}^{\pm}$, correspondingly. This means that, in contrast to the uniform case, there is almost no compensation of partial sums $\rho^{\pm}_n$ in the integral metric on ${\mathbb R}$. The value of $\beta_p$ was revisited in~\cite{[CHDA]} but the question about the sharpness of the constant $\tilde\beta_p$ obtained there with respect to the rate of $p$ still remains open.

Later on, analogues of Gorin's problem were considered for other sets (semi-axes, segments, rectifiable compacts, etc.) and with a different normalisation of~(2). Detailed history of the problems and related results in this direction are summarised in the survey~\cite{[7]}. Here we only mention Gelfond's problem that we consider below. In \cite{[2]}, estimates of type~(1) but with the normalisation of the derivative of SPF are considered,
$$
d'_n(\mathbb
R,p)=\inf\left\{Y(\rho_n):\;\|\rho'_n\|_{L^{p}({\mathbb
        R})}\le 1 \right\},\quad \quad 1<p\le\infty.
 \eqno{(4)}
$$
The following estimate was obtained in \cite{[2]}: $d'_n(\mathbb R,\infty)\ge
{\rm const}\cdot 2^{-n/4}$. Later, Nikolaev generalised and improved some of Gelfond's results. In particular, he showed in \cite{[5]} that 
$$
 d'_n(\mathbb R,\infty) \ge {\rm const}\cdot n^{-3/2}.
$$
Much later, the following weak equivalence was proved in \cite{[6]}:
$$
\inf\left\{Y(\rho_n):\;\|\rho'_n\|_{L^{\infty}({\mathbb
        R})}\le 1,\, \rho_n=\rho_n^{+} \right\}
        \asymp \frac{\ln n}{\sqrt{n}},
 \eqno{(5)}
$$
where, in contrast to~(4), it is additionally assumed that all poles of $\rho_n=\rho_n^{+}$ belong to the half-plane~$\mathbb{C}^{+}$. It is plausible that the same holds in the general case, for $d'_n(\mathbb R,\infty)$.

Note that the exchange $\varrho(z)=c\rho_n(c\,z)$,
$c=\|\rho_n\|^{-q}_{L^{p}({\mathbb R})}$ saves the form of a SPF, and $\|\varrho\|_{L^{p}({\mathbb R})}=1$ and
$Y(\varrho)=c^{-1}Y(\rho_n)$. Consequently, (1) may be rewritten as
$$
d_n(\mathbb R,p)=\inf_{\rho_n}\left\{Y(\rho_n)\,
\|\rho_n\|^{q}_{L^{p}({\mathbb R})} \right\},\qquad
p^{-1}+q^{-1}=1,
  \eqno{(6)}
$$
where the infimum is taken over all SPFs~(2) with no poles on~${\mathbb R}$. Analogously, considering the SPF
$$
\varrho(z)=c\rho_n\left(c\,z)\right),\qquad c:=
 \|\rho'_n\|^{-\frac{q}{q+1}}_{L^{p}({\mathbb R})},
         $$
one gets $\|\varrho'\|_{L^{p}({\mathbb R})}=1$ and
$Y(\varrho)=c^{-1}Y(\rho_n)$, and therefore 
$$
d'_n(\mathbb R,p)=\inf_{\rho_n}\left\{Y(\rho_n)\,
\|\rho_n'\|^{\frac{q}{q+1}}_{L^{p}({\mathbb R})} \right\},\qquad
p^{-1}+q^{-1}=1. 
\eqno{(7)}
$$

Thus Gorin's and Gelfond's problems can be thought as finding the least deviation from zero in $L^{p}({\mathbb R})$ of SPFs~(2) and their derivatives under the condition $Y(\rho_n)=1$, or, which is the same, under the condition that all SPFs~(2) have a common fixed pole, say, ${\xi}_1=i$. In this sense, the problems are analogues of classical Chebyshev's problem on the least deviation from zero of a unitary polynomial of a fixed degree. This circumstance, in particular, leads to more general approximation problems for SPFs (2) and their derivatives on ${\mathbb R}$ and other sets, making the estimates for (1) and (4) still topical (see \cite{[7]}).

Recall that the two-sided estimate (3) is valid for the class of {\it all} SPFs~(2), with no attention to the multiplicity of the roots of $Q$. A natural question about the estimation of $Y(\rho_n)$ for an individual normalised SPF $\rho_n$, taking into account $n_k$, arises. It is answered in the following theorem.

\medskip

{\bf Theorem 1.} {\it There is an absolute $c>0$ such that for any pole ${\xi}_k$ it holds that}
$$
|{\rm Im}\,{\xi}_k|\cdot\|\rho_n\|_{L^{\infty}({\mathbb
        R})}\ge c\; {\frac { \left({\ln n}\right)
        ^{1/n_{k}}+1}{ \left( {\ln n} \right) ^{1/n_{k}}-1}}\cdot \frac{
    \ln {\ln n}}{{\ln n}}>2c\; \frac {n_k}{\ln n}, \qquad n\ge 4.
    \eqno{(8)}
$$

\medskip

Note that the second inequality in (8) follows from the simple inequality
$$
\frac{\mu^{t}+1}{\mu^{t}-1}>\frac{2}{t\ln
\mu},\qquad \text{where}\qquad \mu=\ln n,\quad t=\frac{1}{n_k}>0.
$$

Thus Theorem 1 provides a continuous scale of additional factors in~(3). For example, if $n_k\le \ln\ln n$, then the first inequality in~(8) has the same rate as in (3). For $n_k$, satisfying the opposite inequality, the second inequality in (8) is more precise than~(3).

\smallskip

As for the estimates for $d'_n(\mathbb R,\infty)$, we prove the following theorem in the general case, i.e. without any assumptions on the location of poles.

\medskip

{\bf Theorem 2.} {\it There exists an absolute $c>0$ such that}
$$
d' _n(\mathbb R,\infty)\ge c\,\sqrt{\frac{\ln n}{n}},\qquad n\ge
n_0.
 \eqno{(9)}
$$

\section{Proof of the estimate (8)}

{\bf 2.1. Assumptions.} It is sufficient to prove (8) in the case when $\|\rho_n\|_{L^{\infty}({\mathbb R})}=1$ (see (6)). For determinacy, we obtain a lower estimate for $y_1={\rm Im}\, {z}_1$ assuming that ${z}_1=iy_1$ is one of the poles of $\rho_n$ belonging to the upper half-plane ${\mathbb C}^+$.

First we get estimates under the following additional assumptions.

1) The poles of $\rho_n$ and corresponding residues are symmetric with respect to the real and imaginary axes so that the poles on the imaginary axis have even residues. This happens e.g. if the symmetrisation from Section~2.4 is applied.

By $z_k$, $k=1,\ldots,m$, we denote the poles of SPF $\rho_n$ belonging to the upper half-plane ${\mathbb C}^+$, and by $n_k$ the corresponding residues so that the order of SPF equals
$n=2\sum_{k=1}^m n_k$. We aim to estimate the imaginary part $y_1>0$ of the pole $z_1=iy_1$.

Let
$$
B(z):=
 \prod_{k=1}^{m}\frac{(z-z_k)^{n_k}}{(z-\overline{z_k})^{n_k}},
 \qquad \mu(x)=\frac{1}{2i}\frac{B'(x)}{B(x)}=
 \sum_{k=1}^{m}\frac{n_k\,y_k}{(x-x_k)^2+y_k^2},\quad x\in
{\mathbb{R}}.
 \eqno{(10)}
$$

We use one more assumption.

2) For real $x_1$ and $x_2$
$$
|\mu(x_1)-\mu(x_2)|\le 3 \ln\left(1+\frac{r}{2y_1}\right),\qquad
r:=|x_1-x_2|.
 \eqno{(11)}
$$
Below we show that the assumptions 1)--2) do not limit the generality of the problem.

In the Blaschke product (10), the products of factors with poles $\overline{z_k}$ and $-z_k$, being symmetric with respect to the imaginary axis, are non-negative on the imaginary axis as
$$
\frac{(iy-z_k)}{(iy-\overline{z_k})}
\frac{(iy+\overline{z_k})}{(iy+z_k)}=
\frac{(iy-z_k)}{(iy-\overline{z_k})}\overline{
\frac{(-iy+{z_k})}{(-iy+\overline{z_k})}}=
\left|\frac{iy-z_k}{iy-\overline{z_k}}\right|^2.
$$
Therefore the symmetry assumption 1) implies that $B(iy)\ge 0$ for all $y\in
{\mathbb R}$ and $0<B(iy)<1$ for $y>0$. The assumption 1) also implies that  $\mu$ is an even positive function on the real axis. Furthermore, the following partial fraction decomposition holds:
$$
\frac{1-B(z)}{1+B(z)}=i\sum_{k=1}^{2\eta}\frac{1}{\mu
(t_{k})}\frac{1}{z-t_{k}},\qquad 2\eta:=n/2=\sum_{k=1}^m n_k,
 \eqno{(12)}
$$
with pairwise distinct finite real $t_k$ being the roots of the equation $B(x)=-1$. The points $t_k$ locate on the real axis symmetrically with respect to the origin (it follows from the equality $B(x)=\overline{B(-x)}$). For determinacy, let $t_k<t_{k+1}$ ($k=1,\ldots,2\eta-1$) and let $t_k<0$ for $k=1,\ldots,\eta$ and $t_k>0$ for $k=\eta+1,\ldots,2\eta$. Set
$$
r_k=t_{\eta+k},\qquad k=1,\ldots,\eta.
$$
On each segment $[t_k,t_{k+1}]$ the argument of the Blaschke product $B(x)$ has increment of  $2\pi$, in particular,
$$
\int_0^{r_k}\mu(x)\,dx=\frac{1}{2i}\int_0^{r_k}(\ln(B(x)))'\,dx
=\frac{1}{2}\int_0^{r_k}(\arg(B(x)))'\,dx= \frac{\pi}{2}(2k-1).
 \eqno{(13)}
$$
Fix $\theta\in (0,1)$ and $y_0=y_1\theta $ ($y_1={\rm
Im}\,z_1$). Since
$$
0<B(iy_0)=\frac{(1-\theta)^{n_1}}{(1+\theta)^{n_1}}
 \prod_{k=2}^{m}\frac{(iy_0-z_k)^{n_k}}{(iy_0-\overline{z_k})^{n_k}}
<\varepsilon,\qquad
\varepsilon:=\frac{(1-\theta)^{n_1}}{(1+\theta)^{n_1}},
 \eqno{(14)}
$$
the decomposition (12) leads to
$$
\frac{1-\varepsilon}{1+\varepsilon}
 <\frac{1-B(iy_0)}{1+B(iy_0)}=i\sum_{k=1}^{2\eta}\frac{1}{\mu
(t_{k})}\frac{1}{iy_0-t_k}=\sum_{k=1}^{\eta}\frac{2}{\mu
(r_{k})}\frac{y_0}{y_0^2+r_k^2},
$$
and thus 
$$
\sum_{k=1}^{\eta}\frac{2}{\mu (r_{k})}\frac{y_0}{y_0^2+r_k^2}\ge
1-\delta,\qquad \delta:=\frac{2\varepsilon}{1+\varepsilon}
 =\frac{2(1-\theta)^{n_1}}{(1-\theta)^{n_1}+(1+\theta)^{n_1}}.
 \eqno{(15)}
$$

\smallskip

{\bf 2.2. Estimate for the sum in (15)}. Let $\tau>0$ and
$$
\mu_1(\tau)=\min_{[0,\tau]}\mu(x),\qquad
\mu_2(\tau)=\max_{[0,\tau]}\mu(x).
$$
Fix $r>0$ and divide the sum in (15) into the two:
$$
S_1(r)+S_2(r):=\left(\sum_{r_k\le
r}+\sum_{r_k>r}\right)\frac{2}{\mu
(r_{k})}\frac{y_0}{y_0^2+r_k^2}.
$$
To estimate $S_1$, take into account (13):
$$
r_k\mu_2(r_k)\ge \int_0^{r_k}\mu(x)\,dx=\frac{\pi}{2}(2k-1),\qquad
r_k\ge \frac{{\pi}(2k-1)}{2\mu_2(r_k)},
$$
which implies that
\begin{align*}
S_1(r)&\le \sum_{r_k\le r}\frac{2}{\mu_1 	(r)}\frac{y_0}{y_0^2+\frac{{\pi^2}(2k-1)^2}{4\mu_2^2(r)}}\\
	& \le\frac{4\mu_2(r)}{\mu_1
	(r)}\sum_{k=1}^{\infty}\frac{2\mu_2(r)y_0}{\left(2\mu_2(r)y_0\right)^2+{{\pi^2}(2k-1)^2}}\\
&=\frac{\mu_2(r)}{\mu_1(r)}\frac{e^{2\mu_2(r)y_0}-1}{e^{2\mu_2(r)y_0}+1}.
\end{align*}
The sum of the series is known, see e.g. \cite[Chapter 2, \S 3]{Fiht}.

To estimate $S_2$, consider (12). By Cauchy's integral formula,
$$
\frac{1}{2\pi}\int_{-i\infty}^{i\infty}\frac{1}{(\xi+r)^2}
\frac{1-B(\xi)}{1+B(\xi)}\,d\xi=
\sum_{k=1}^{\eta}\frac{1}{\mu(r_k)}\frac{1}{(r+r_k)^2}.
$$
Furthermore, recall that $B(iy)\ge 0$ for all $y\in {\mathbb R}$. By this reason,
the modulus of the left hand side of (12) is at most $1$ and therefore
$$
\frac{1}{2\pi}\left|\int_{-i\infty}^{i\infty}\frac{1}{(\xi+r)^2}
\frac{1-B(\xi)}{1+B(\xi)}\right|\le
\frac{1}{2\pi}\int_{-i\infty}^{i\infty}\frac{1}{|\xi+r|^2}=
\frac{1}{2r}.
$$
This implies that
$$ S_2(r)\le 2\sum_{r_k>r}\frac{y_0}{\mu (r_{k})r_k^2}\le
8\sum_{r_k>r}\frac{1}{\mu (r_{k})}\frac{y_0}{(r+r_k)^2}\le
\frac{4y_0}{r}.
$$

The sum of the estimates for $S_1$ and $S_2$ and (15) give the inequity
$$
\frac{\mu_2(r)}{\mu_1(r)}\frac{e^{2\mu_2(r)y_0}-1}{e^{2\mu_2(r)y_0}+1}
+\frac{4y_0}{r}\ge 1-\delta,\qquad y_0=y_1\theta.
$$
Thus the following lemma is true. 
\medskip

{\bf Lemma 1.} {\it Under the assumptions $1)$ and $2)$, for any $r>0$ it holds that 
$$
e^{2\theta\,\mu_2(r)y_1}\ge
\frac{\mu_2(r)+\mu_1(r)-\delta\mu_1(r)-4y_0\mu_1(r)/r}
 {\mu_2(r)-\mu_1(r)+\delta\mu_1(r)+4y_0\mu_1(r)/r}, \eqno{(16)}
$$
where
$$
r>0,\quad \theta\in (0,1),
 \quad \delta=\frac{2(1-\theta)^{n_1}}{(1-\theta)^{n_1}
 +(1+\theta)^{n_1}}.
$$
}

\smallskip
{\bf 2.3. Choice of $r$ and $\delta\in (0,1)$}. Since
$\mu_2(r)\ge \mu_1(r)\ge 0$, (16) gives
$$
e^{2\theta\,\mu_2(r)y_1}\ge \frac{\mu_2(r)-4y_0\mu_2(r)/r}
 {\mu_2(r)-\mu_1(r)+\delta\mu_2(r)+4y_0\mu_2(r)/r}.
$$

From now on we think that $y_1\le n_1/10$ (in the otherwise case the inequality (8) is obvious), therefore $\mu_2(r)>\mu(0)\ge{n_1}/{y_1}>
10$. Choose $r>0$ and then $\delta$ from the conditions
$$
r=4\mu_2(r)y_1,\qquad \mu_2(r)\delta=1.
$$
Such $r$s obviously exist, possibly they are multiple. Solving the latter equation with respect to $\theta$ gives
 $$
 \theta=\frac{(2\mu_2(r)-1)^{1/n_1}-1}{(2\mu_2(r)-1)^{1/n_1}+1}.
  \eqno{(17)}
 $$
For the chosen $r$ and $\theta$ the following inequalities are valid:
$$
e^{2\theta\,\mu_2(r)y_1}\ge \frac{\mu_2(r)-1}
 {\mu_2(r)-\mu_1(r)+2},\qquad
y_1\ge \frac{1}{2\theta\,\mu_2(r)}\ln\left(\frac{\mu_2(r)-1}
 {\mu_2(r)-\mu_1(r)+2}\right),
$$
$$
y_1\ge \frac{1}{2\,\mu_2(r)}
\frac{(2\mu_2(r)-1)^{1/n_1}+1}{(2\mu_2(r)-1)^{1/n_1}-1}
 \ln\left(\frac{\mu_2(r)-1}
 {\mu_2(r)-\mu_1(r)+2}\right).
$$
By (11), the choice of $r=4\mu_2(r)y_1$ and the inequality $\mu_2(r)>10$, 
$$
\mu_2(r)-\mu_1(r)\le 3\ln\left(1+\frac{r}{2y_1}\right)\le
3\ln\left(1+2\mu_2(r)\right)< 4\ln \mu_2(r).
$$

Thus we have proved the following lemma.

\medskip
{\bf Lemma 2.} {\it If $r=4\mu_2(r)y_1$ and $y_1\le n_1/10$, then}
$$
y_1\ge \frac{1}{2\,\mu_2(r)}
\frac{(2\mu_2(r)-1)^{1/n_1}+1}{(2\mu_2(r)-1)^{1/n_1}-1}
\ln\left(\frac{\mu_2(r)-1}
{2+4\ln\mu_2(r)}\right)
$$
$$
\ge c\,
\frac{\mu_2(r)^{1/n_1}+1}{\mu_2(r)^{1/n_1}-1}\cdot
\frac{\ln \mu_2(r)}{\mu_2(r)}.\eqno{(18)}
$$

Now recall that $\|\rho_n\|_{L^{\infty}({\mathbb R})}\le
1$ and use the following estimate from \cite{[8]}:
$$
\|\mu\|_{L^{\infty}({\mathbb R})}\le {\rm const}\,\ln n.
$$
Taking into account that the minorant in (18) is decreasing as a function of $\mu_2(r)$ (it can be easily checked) and that $\mu_2(r)\le \|\mu\|_{L^{\infty}({\mathbb
R})}$, we come to the first inequality in (8). Thus the inequality (8) is proved under the assumptions 1)--2).

\smallskip
{\bf 2.4. General case}. The general case of SPF (2) (of a given order $n$) can be reduced to Lemma~2 as follows. Considering $z_1=iy_1$ and
$y_1>0$, we make two symmetrisations of the form
$$
s_1(z)=\rho_n(z)+\overline{\rho_n(\bar z)},\quad \sigma_{0}(z)=
s_1(z)-\overline{s_1(-\bar z)}.
$$
These give a SPF of the form $\sigma_{0}(z)=\sigma(z)+\overline{\sigma(\bar
z)}$ of order $4n$ with poles symmetric with respect to the real and imaginary axes so that the poles of $\sigma$ belong to ${\mathbb
C}^+$ and the poles of $\overline{\sigma(\bar
z)}$ to ${\mathbb C}^{-}$. Obviously, the maximum of $|\sigma_{0}(x)|$ is at most four times more that the maximum of $|\rho_n(x)|$ and therefore 
$\|\sigma_{0}\|_{L^{\infty}({\mathbb R})}\le 4$. One of the poles of $\sigma_{0}$ is still $z_1=iy_1$, $y_1>0$, with the residue $\ge 2n_1$. Furthermore, we exchange  $\sigma_{0}(z)$ for SPF $\rho(z)=\sigma(z-iy_1)+\overline{\sigma(\overline{z+iy_1)}}$, i.e. move the poles of $\sigma_{0}(z)$ from the real axis by $y_1$ so that one of the poles of $\rho$ is
$2z_1=2iy_1$, with the residue $\ge 2n_1$. The value of $\|\rho\|_{L^{\infty}(\mathbb R)}$ is then at most four times more than the maximum modulus of the initial SPF (2). This follows from the maximum modulus principle for subharmonic functions:
$$
|\rho(x)|=2|{\rm Re}\,\sigma(x-iy_1)|\le 2\|{\rm Re}\,\sigma\|_{L^\infty(\mathbb
R)}= \|\sigma_{0}\|_{L^\infty(\mathbb R)}\le 4.
$$

Now let us show that the SPF $R$ defined by
$$
R(z):=4^{-1}\rho(4^{-1}z)=
 \varrho(z)+\overline{\varrho(\overline{z})},\qquad
 \varrho(z):=\frac{1}{4}\sigma\left(\frac{z-iy_1}{4}\right),
  \eqno{(19)}
 $$
satisfies the assumptions 1)--2). Indeed, its $\sup$-norm on
$\mathbb R$ is at most $1$, its poles and residues are symmetric with respect to the coordinate axes and moreover the residue of its pole $8z_1=8y_1i$ is at least  $2n_1$. What is more,  Cauchy's integral formula for $z\in {\mathbb C}^-$ gives
$$
|\sigma'(z)|\le \frac{1}{2\pi}\int_{\mathbb
R}\frac{|\sigma_{0}(x)|}{|x-z|^2}\,dx\le \frac{4}{|{\rm Im}\, z|}
$$
and therefore by (19),
$$
|\varrho'(z)|\le \frac{1}{y_1+|{\rm Im}\, z|}.
$$
By integrating this estimate against the right angle $\gamma$ with equal sides that belongs to the lower half-plane and is based on the segment $[x_1,x_2]$, $r=x_2-x_1>0$, we get
$$
|\varrho(x_1)-\varrho(x_2)|\le \int_{\gamma}\frac{|dz|}{y_1+|{\rm
Im}\, z|}
  =2\sqrt{2}\int_{0}^{r/2}\frac{dx}{-x+y_1+r/2}
  <3\ln\left(1+\frac{r}{2y_1}\right).
$$

Thus the inequality (11) holds and the SPF $R$ satisfies the assumptions 1)--2). Consequently, the estimate (8) is valid for it. The difference between $R$ and the initial SPF $\rho_n$ by means of $n$, $n_k$ and $\xi_k$ does not influence the rate in the estimate~(8) but only changes the absolute constant $c$.

\section{Proof of the estimate (9)}

Let $\rho_n(z)=\rho^{+}(z)+\rho^{-}(z)$, where $\rho^{\pm}$ are SPFs whose poles lie in ${\mathbb C}^{\pm}$. Let
$$
\sigma(z)=\rho'_n(z),\quad \sigma_1(z)=(\rho^{+}(z))',\quad
\sigma_{2}(z)=(\rho^{-}(z))',
$$
so that $ \sigma(z)=\sigma_1(z)+\sigma_2(z)$. For simplicity,  suppose that $\|\sigma\|_{L^{\infty}({\mathbb
R})}=1$.

\smallskip

{\bf Lemma 3}. {\it  Given a fixed $n\ge 2$ and $\|\sigma\|_{L^{\infty}({\mathbb R})}=1$,}
$$
\|\sigma_1(\cdot-ih)\|_{L^{\infty}({\mathbb R})}\le 5\ln n,\qquad
h=\frac{1}{n^2}.
  \eqno{(20)}
$$
{\bf Proof.} Cauchy's integral formula for $z\in
{\mathbb C}^-$ gives
$$
|\sigma_1'(z)|\le \frac{1}{2\pi}\int_{\mathbb
R}\frac{|\sigma(x)|}{|x-z|^2}\,dx\le
\frac{\|\sigma\|_{L^{\infty}({\mathbb R})}}{|{\rm Im}\,
z|}=\frac{1}{|{\rm Im}\, z|}.
$$
For a fixed $x\in\mathbb R$ and $h=1/n^2$, this implies that
$$
|\sigma_1(x-ih)-\sigma_1(x-i/h)|\le
\int_{h}^{1/h}\frac{dy}{y}=4\ln n.
$$
Obviously, $|\sigma_1(x-i/h)|\le n^{-3}$ and therefore
$$
|\sigma_1(x-ih)|\le n^{-3}+4\ln n<5\ln n, \qquad n\ge 2.
$$
The inequality (20) is proved.

\medskip

Furthermore, the estimate (5) and definition (7) for $p=\infty$, $q=1$, leads to
$$
(h+Y(\rho^{+}))\,
\|\sigma_1(\cdot+ih)\|^{\frac{1}{2}}_{L^{\infty}({\mathbb R})}\ge
c_1\frac{\ln n}{\sqrt{n}},
$$
so that Lemma~3 implies
$$
h+Y(\rho^{+}) \ge \frac{c_1}{\sqrt{5}}\frac{\sqrt{\ln
n}}{\sqrt{n}},\qquad Y(\rho^{+}) \ge
\frac{c_1}{\sqrt{5}}\frac{\sqrt{\ln
n}}{\sqrt{n}}-\frac{1}{n^2}>c\frac{\sqrt{\ln n}}{\sqrt{n}},\quad
n\ge n_0(c_1).
$$
Analogous inequalities hold for $Y(\rho^{-})$, too, thus the required inequality (9) follows.

\section{Acknowledgements}

The reported study was funded by Russian Ministry of Education and Science
(task number 1.574.2016/1.4) and RFBR (project number 18-01-00744).


\end{document}